\documentclass[11pt, reqno] {amsart}
\usepackage {fullpage, amsmath, amssymb,latexsym,amsfonts,amsthm,mathtools,ulem, graphics, setspace}
\usepackage{csquotes}
\newtheorem{theorem}{Theorem}

\newtheorem{remark}{Remark}
\newtheorem{proposition}{Proposition}

\newtheorem{problem}{Problem}
\newtheorem{assume}{Assumption}
\DeclareMathOperator{\re}{Re}

\usepackage[colorlinks=true, pdfstartview=FitV, linkcolor=blue,citecolor=red, urlcolor=blue]{hyperref} 

\usepackage{graphicx}
\usepackage{caption}
\usepackage{subcaption}
\usepackage{float}

\usepackage{pdfpages}

\newcommand{\mathsym}[1]{{}}
\newcommand{\unicode}[1]{{}}

\usepackage{listings}
\begin {document} 
\title{Short mollifiers of the Riemann zeta-function}  

\author{J.~B. Conrey}
\address{American Institute of Mathematics \\ 1200 E California Blvd, Pasadena CA 91125}
\email{\href{mailto:conrey@aimath.org}{conrey@aimath.org}}

\author{D.~W. Farmer}
\address{American Institute of Mathematics \\ 1200 E California Blvd, Pasadena CA 91125}
\email{\href{mailto:farmer@aimath.org}{farmer@aimath.org}}

\author{C.-H. Kwan}
\address{University College London \\ 25 Gordon Street, London WC1H 0AY, United Kingdom }
\email{\href{mailto:ucahckw@ucl.ac.uk}{ucahckw@ucl.ac.uk}}

\author{Y. Lin}
\address{Data Science Institute, Shandong University, Jinan 250100, China}
\email{\href{mailto:yongxiao.lin@sdu.edu.cn}{yongxiao.lin@sdu.edu.cn}}

\author{C.~L. Turnage-Butterbaugh}
\address{Carleton College \\ 1 North College Street Northfield, MN 57707, USA}
\email{\href{mailto:cturnageb@carleton.edu}{cturnageb@carleton.edu}}



\date{\today}









\begin{abstract}
We apply the calculus of variations to construct a new sequence of linear combinations of derivatives of the Riemann $\zeta$-function adapted to Levinson’s method, which yield a positive proportion of zeros of the $\zeta$-function on the critical line, regardless of how short the mollifier is. Our construction extends readily to modular $L$-functions.  Even with Levinson's original choice of mollifier, our method more than doubles the proportions of zeros on the critical line for modular $L$-functions previously obtained by Bernard and K\"uhn--Robles--Zeindler, while relying on the same arithmetic inputs.  This indicates that optimizing the linear combinations, an approach that has received relatively little attention, has a more pronounced effect than refining the mollifier when it is short. Curiously, our linear combinations provide non-trivial smooth approximations of Siegel’s $\mathfrak{f}$-function in the celebrated Riemann--Siegel formula.  
\end{abstract}

\maketitle

\section{Introduction}\label{sect: intro}
There are two well-known methods for proving that a positive proportion of the zeros of the Riemann zeta-function $\zeta(s)$ ($s=\sigma+it$) lie on the critical line $\sigma =1/2$: Selberg's method \cite{Sel42} and Levinson's method \cite{Lev74}. Both of these methods involve the use of {\it mollifiers}.  In Selberg's method, roughly speaking,  one calculates a mollified second moment of $\zeta(s)$ {\it on} the critical line, where the mollifier itself is a fourth power of an approximation to $\zeta(s)^{1/2}$; see \cite[Chapter 24.2]{IK04}.  In Levinson's method, one calculates a second mollified moment of $\zeta(s)$ {\it off} the critical line, where the mollifier is an approximation to $\zeta(s)^{-1}$.

\begin{figure}
    \includegraphics[scale=0.35]{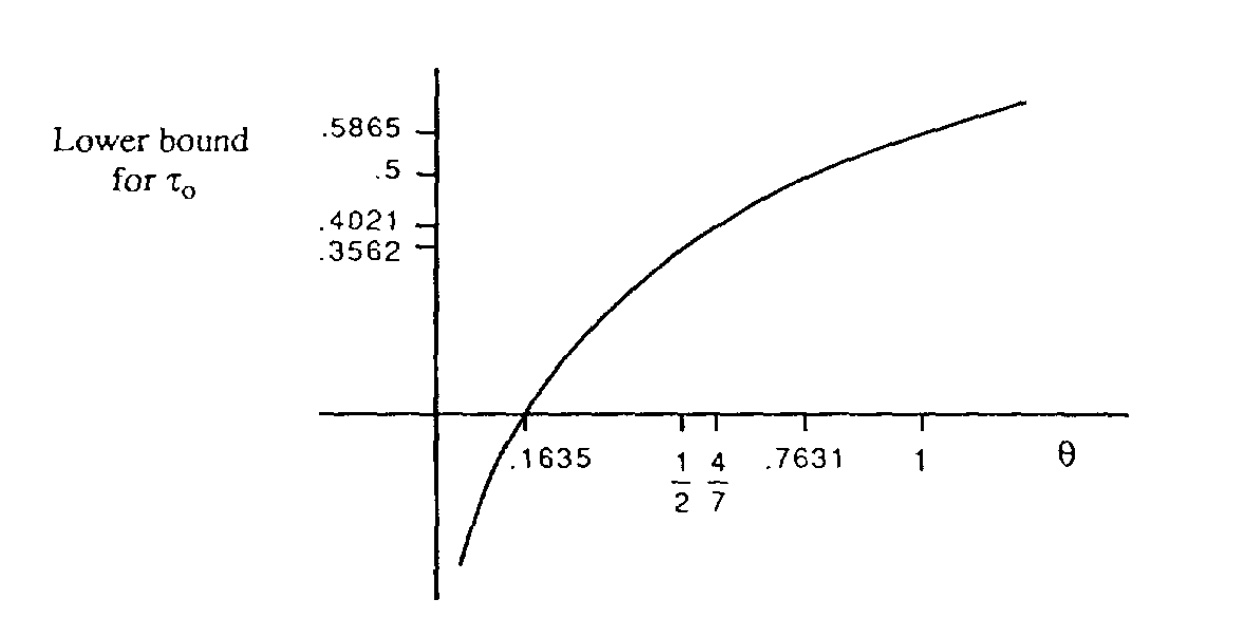}
    \caption{\footnotesize Proportion of zeros on the critical line as a function of the length $\theta$ via Levinson’s method; reprinted from \cite[p. 215]{Far94}}
    \label{fig: farmer}
\end{figure}

Selberg's method has the advantage of working regardless of the length of the mollifier. By contrast, it has long been believed that Levinson's method suffers from a major drawback:  to obtain a positive proportion of zeros on the critical line, the mollifier must be sufficiently long. This viewpoint is supported Farmer  \cite[p. 215]{Far94}, who provided the plot shown in Figure \ref{fig: farmer}
\noindent and remarked,  ``It is interesting that Levinson's method does not yield a result if $\theta$ is too small.'' Levinson's method, however, tends to yield much stronger explicit bounds for the proportion, and in some cases, can also establish a positive proportion of simple zeros on the critical line. As described in \cite[Chapter 16]{Iwa14}, \begin{displayquote}
Levinson’s approach seems to be a gamble, because it is not clear up front that at the end the numerical constants are good enough to yield a positive lower bound for [the number of critical zeros up to height $T$].
\end{displayquote}

In this paper, we prove that Levinson's method, as modified by Conrey \cite{Con89}, will in fact give a positive proportion of zeros on the critical line no matter how short the mollifier is---contrary to the common belief.  

\subsection{Levinson's method}\label{introLev}
Levinson's starting point is a theorem of Speiser \cite{Spe35}, which states that the Riemann Hypothesis is equivalent to the assertion that all non-real zeros of $\zeta'(s)$ are in the half-plane $\sigma\ge 1/2$. This result was refined by Levinson and Montgomery \cite{LM74}, who proved a (essentially) one-to-one correspondence between the non-real zeros of $\zeta(s)$ and those of $\zeta'(s)$ in the region $\sigma<1/2$ (presumably both sets are empty). 

Subsequently,  Levinson \cite{Lev74} proved that the number of zeros of $\zeta'(s)$ in the region $\sigma<1/2$, $0< t<T$ is at most $(1/3)(T/2\pi)\log T$, i.e. at most $1/3$ of the zeros are to the \textit{left} of $\sigma=1/2$. The same result holds for $\zeta(s)$ by the aforementioned theorem of \cite{LM74}. The functional equation of $\zeta(s)$ allows us to deduce that at most $1/3$ of the zeros of $\zeta(s)$ are to the \textit{right} of $\sigma=1/2$.  This leaves \textit{at least} $1/3$ of the zeros of $\zeta(s)$  \textit{on} the line $\sigma=1/2$. 

Using the functional equation again, the zeros of $\zeta'(s)$ to the left of $\sigma=1/2$ correspond to the zeros of $\zeta(s)+L^{-1}\zeta'(s)$ ($L:= \log (T/2\pi)$) on the \textit{right} of $\sigma=1/2$. This formulation turns out to be more practical. The essential ingredient to Levinson's argument is an asymptotic formula for the moment
\begin{align}\label{molsecmofir}
    \frac{1}{T} \, \int_0^T |(\zeta+L^{-1}\zeta')(a+it)|^2 |M(a+it)|^2 ~dt,
\end{align}
where 
$a$ is a real number slightly less than $1/2$ (see (\ref{eqn: aleft1/2log})), and $M$ is a \textit{mollifier} of the form
\begin{align}\label{eqn: LevmollP=x}
    M(s) \, = \, \sum_{n\le y}\frac{\mu(n)}{n^{s+1/2-a}} \frac{\log (y/n)}{\log y},
\end{align}
where $\mu(n)$ is the M\"obius $\mu$-function, $y=T^\theta$,
and $\theta>0$ is the `length' of the mollifier, the key parameter. 

Levinson established (\ref{molsecmofir}) for all $\theta<1/2$. The admissible range for $\theta$ was later extended by Conrey \cite{Con89}  to $\theta < 4/7$.  With refinements to be discussed below, he showed that at least $40\%$ of the zeros of $\zeta(s)$ are on the critical line \cite[Theorem 1]{Con89}. In \cite{Far93}, Farmer proposed the `$\theta=\infty$' conjecture, that is, the admissible value for $\theta$ can be arbitrarily large. He showed that this would imply $100\%$ of the zeros of $\zeta(s)$ lie on $\sigma=1/2$. More recently,  Bettin and Gonek \cite{BG17} proved that `$\theta=\infty$' in fact implies the Riemann Hypothesis.

Conrey \cite{Con89} refined Levinson’s method in two crucial ways. First, he observed that the role of $\zeta(s)+L^{-1}\zeta'(s)$ in (\ref{molsecmofir}) can be replaced by any \textit{linear combination} of derivatives of $\zeta(s)$ of the form $Q(-\frac{1}{L}\frac{d}{ds})\zeta(s)$, where $Q$ is a real polynomial that satisfies 
\begin{align}\label{functeqn}
    Q(0) \, = \, 1 \hspace{15pt} \text{and} \hspace{15pt} Q'(y) \, =\,  Q'(1-y). 
\end{align}
The condition (\ref{functeqn}) stems from the functional equation of $\zeta(s)$. Second, he considered a more general class of mollifiers of the form
\begin{align}
    M(s,P):=\sum_{n\le y}\frac{\mu(n)}{n^{s+1/2-a}} P\Big( \frac{\log (y/n)}{\log y}\Big),
\end{align}
where $P$ is a real polynomial with 
\begin{align}\label{molliP}
    P(0) \, = \, 0 \hspace{15pt}  \text{and}  \hspace{15pt} P(1) \, = \, 1.
\end{align}
 He established an asymptotic formula for the corresponding mollified moment (see (\ref{Conmollif}) and (\ref{eqn cpqroptim})) for any such $P$, $Q$ and $\theta<4/7$, and determined the optimal choice of $P$.


\subsection{Results}

We prove that no matter how small $\theta>0$ is in (\ref{eqn: LevmollP=x}), there exists an admissible linear combination
of derivatives of $\zeta(s)$ that yields a positive proportion of zeros of $\zeta(s)$ on the critical line via Levinson's method. Moreover, this proportion can be made quantitative; see \textbf{Theorem \ref{precimain}}. In fact, our linear combination is optimal with respect to the variational problem we consider (see Problem \ref{mainprob} and equation (\ref{eqn:S})).  Such an optimizer (for $Q$) is significantly more intricate than the one for the mollifier polynomial $P$ in \cite[p. 9]{Con89}.

For small values of $\theta$, the improvement gained by choosing $Q$ strategically is significant. For example, in Bernard's work \cite{Ber15} on critical zeros of modular (degree two) $L$-functions, he obtains a proportion of 2.97\% of zeros on $\sigma=1/2$ (using the optimal $P$ and an empirical $Q$ of degree $7$), whereas we get 6.32\% (using $P(x)=x$ and our newly constructed $Q$), both unconditionally and using the same shifted convolution estimate. We have not attempted to optimize the proportion here (a method for full optimization is described in Section \ref{sect: gene}), but our result already comes close to the current record achieved in \cite{AT21} (using the optimal $P$ and an empirical $Q$ of degree $27$!). 

The broader goal is to initiate a systematic and \textit{theoretical} study of how to select effective linear combinations of derivatives of the $\zeta$- (or $L$-)functions in analytic applications, in contrast to the predominantly \textit{empirical}, computer-experimental approach used to date.


\subsection{Connections to Siegel's method}\label{sect: Siegmethod}

Let $h(s):= \pi^{-s/2}\Gamma(s/2)$. Levinson's method detects zeros of $\zeta(s)$ on $\sigma=1/2$ by identifying when $\mathrm{arg}\,  (h\zeta')(1/2+it) \equiv \pi/2 \, (\bmod\, \pi)$ occurs. By the argument principle, this task reduces to bounding the number of zeros of $\zeta'(s)$ on $\sigma <1/2$. In the introduction to \cite{Lev74}, Levinson noted that the idea of detecting zeros of $\zeta(s)$ on the critical line via changes in argument and moment estimates (of a suitable function) had been used in Siegel’s 1932 work \cite{Sie32}. 



In \cite{Sie32}, instead of working with $\zeta'(s)$ as in \cite{Lev74}, Siegel introduced an interesting entire function:
\begin{align}
    \mathfrak{f}(s) \, := \, \int_{L} \frac{e^{i\pi w^2}}{e^{i\pi w}-e^{-i\pi w}} w^{-s} \, dw,
\end{align}
where the line of integration $L$ has slope $1$, with $\mathrm{Im} \, w$ decreasing, and intersects the real axis at a point between $0$ and $1$. The function $\mathfrak{f}(s)$ originates from the \textit{Riemann--Siegel formula}, which is an integral representation of $\zeta(s)$ given by
\begin{align}\label{RieSieg-intrep}
    h(s)\zeta(s) \, = \, 2\, \re h(s)\mathfrak{f}(s) \hspace{20pt} (\re s \, = \, 1/2). 
\end{align}
This formula first appeared in Riemann’s 1859 memoir \cite{Rie} but remained largely unnoticed in the literature until its rediscovery by Siegel. 

Whenever $\mathrm{arg}\,  (h\mathfrak{f})(1/2+it) \equiv \pi/2 \, (\bmod\, \pi)$ occurs, it follows from (\ref{RieSieg-intrep}) that there is a zero of $\zeta(s)$ on the critical line. In fact, Siegel proved that
\begin{align}\label{Sieg: lowbdd}
    N_{0}(T) \, > \, 2 \, \cdot\,  \#\{s= \sigma+it: \, \mathfrak{f}(s)=0, \, \sigma<1/2, \, 0<t<T\},
\end{align}
where $N_{0}(T)$ is the number of zeros of $\zeta(s)$ with $\sigma=1/2$ with $0<t\le T$. A \textit{lower bound} to the right-hand side of (\ref{Sieg: lowbdd}) can be obtained by evaluating the moment
\begin{align}
  \frac{1}{T}  \int_{1}^{T} \, |\mathfrak{f}(\sigma+it)|^2 \, dt,
\end{align}
which turns out to be $\sim C_{\sigma} \cdot T^{1/2-\sigma}$ for some explicit $C_{\sigma}>0$, and for $\sigma<1/2$ to be chosen.  As a consequence, Siegel established, for $T\gg 1$, that
\begin{align}
    N_{0}(T) \, > \, \frac{3}{8\pi} e^{-3/2} T,
\end{align}
which improves upon the result of Hardy and Littlewood \cite{HL21} that $N_{0}(T)\gg T$. 

The function $\mathfrak{f}(s)$ can be \textit{formally} interpreted in Conrey--Levinson's set-up. Consider 
\begin{align}\label{eqn: Siegelstep}
   Q_{\infty}(y) \, := \,  \begin{cases}
        \ \   1 \hspace{15pt} \text{if} \hspace{10pt} 0 \, \le \,  y \, < \,  1/2\\
       \, 1/2 \hspace{10pt} \text{if} \hspace{20pt} y \, = \, 1/2\\
       \ \  0 \hspace{15pt} \text{if} \hspace{10pt} 1/2 \, < \, y \, \le \, 1.
    \end{cases}
\end{align}    
Let $s=1/2+iT$ with $T\gg 1$, and $L:= \log (T/2\pi)$. By the saddle point method, Siegel deduced, up to a negligible error of  $O(T^{-1/4})$, that
\begin{align}
    \mathfrak{f}(s) \, &= \, \sum_{n\le \sqrt{T/2\pi}} \, \frac{1}{n^s}  = \, \sum_{n=1}^{\infty} \, n^{-s}Q_{\infty}\Big(\frac{\log n}{L}\Big)  \, = \, \sum_{n=1}^{\infty} Q_{\infty}\Big(-\frac{1}{L}\frac{d}{ds}\Big)n^{-s} \, = \, Q_{\infty}\Big(-\frac{1}{L}\frac{d}{ds}\Big)\zeta(s). \nonumber
\end{align}

The step function $Q_{\infty}(y)$ satisfies the condition (\ref{functeqn}) apart from failing to be differentiable at $y=1/2$, and therefore does not qualify as an admissible function in Conrey--Levinson's method. In fact, it was commented in \cite[p. 384]{Lev74} that
\begin{displayquote}
It appears to me that the function $f_{1}(s)$ [i.e., $\mathfrak{f}(s)$] is not amenable to improvement with a mollifier.
\end{displayquote}
As Siegel noted at the end of his paper, his method itself \textit{cannot} yield a bound $N_{0}(T)\gg T(\log\log T)^{1+\epsilon}$. However, we will show that

\begin{proposition}\label{appSige}
For any $\theta>0$, let $Q_{\theta}: [0,1]\rightarrow \mathbb{R}$ be the function constructed in Theorem \ref{precimain}, which is used to define linear combinations of derivatives of $\zeta(s)$. Then as $\theta \to 0+$, the pointwise limit of $Q_{\theta}(y)$ converges to the step function $Q_{\infty}(y)$ for any $y\in [0,1]$. 
\end{proposition}

In other words, behind \textbf{Theorem \ref{precimain}} is the idea of mollifying linear combinations of derivatives of $\zeta(s)$ that \textit{approximate} Siegel's function $\mathfrak{f}(s)$. Proposition \ref{appSige} shows that we are actually using polynomials $Q$ of \textit{arbitrarily large} degree in showing Levinson's method works for mollifiers of any length $\theta>0$. Unlike previous works on mollifiers and critical zeros, a key feature of our approach is that our linear combinations vary with the values of $\theta$.



    
\section*{Acknowledgments}
This work was initiated at the AIM workshop \textit{Delta Symbols and the Subconvexity Problem}, held online during the pandemic from 2-6 November 2020. It is a pleasure to thank the AIM staff and organizers for their excellent support. The research was also supported by NSF DMS-1854398 FRG. The third author was supported by the EPSRC grant: EP/W009838/1. The fourth author was supported by the National Key R\&D Program of China (No. 2021\allowbreak YFA1000700). The fifth author is partially supported by NSF CAREER DMS-2239681.

All computations and plots in this article were carried out using \textit{Mathematica}. The relevant commands are included in Appendix \ref{sect: mathematica}. The corresponding dataset (.mx files) and notebooks (.nb files) are available at \url{https://github.com/davidfarmer/shortmollifiers}.

     
\section{Preliminary work}\label{sec: prelim}

\subsection{Set up}

Let $P$ be a polynomial satisfying (\ref{molliP}), $Q$ a polynomial satisfying (\ref{functeqn}), and let $R>0$ be any constant. According to \cite[eq. (39)]{Con89}, if $\theta>0$ is an admissible constant such that the asymptotic formula
\begin{align}\label{Conmollif}
    \frac{1}{T} \int_{0}^{T} \, \Big|Q\Big(-\frac{1}{L} \frac{d}{ds}\Big)\zeta(a+it) M(a+it, P)\Big|^2 \, dt \, \sim \,  c(P,Q,R)
\end{align}
holds as $T\to \infty$, where 
\begin{align}\label{eqn: aleft1/2log}
    a \, = \, 1/2-R/\log (T/2\pi),
\end{align}
and $c(P,Q,R)$ is a constant depending on $P, Q,$ and $R$, then the proportion of zeros of $\zeta(s)$ on the critical line is \textit{at least}  
\begin{align}\label{kappaprop}
    \kappa \, := \, 1-\frac{\log c(P,Q,R)}{R}.
\end{align}
By \cite[Theorem 2]{Con89}, the constant $c(P,Q,R)$ is given by
\begin{eqnarray}\label{eqn cpqroptim}
 1+\frac1 \theta \int_0^1\int_0^1 \big(w(y)P'(x)+\theta w'(y) P(x)\big)^2 ~dx ~dy, 
\end{eqnarray}
where
\begin{equation}\label{eqn cpqroptimwy}
w(y) \, := \, e^{Ry}Q(y).
\end{equation}


To simplify the exposition and illustrate the strength of our construction, we adopt Levinson’s original choice:
\begin{align}\label{linearP}
    P(x) \, = \,  x.
\end{align}
Also, observe that the class of admissible functions $Q$ can be extended to include all continuously differentiable ($C^1$) functions $Q: [0, 1]\rightarrow \mathbb{R}$ with
   \begin{align}\label{Q-funcopt}
        Q(0)=1 \hspace{15pt} \text{and} \hspace{15pt} Q(y) +Q(1-y) \, = \, 1.
   \end{align}
Indeed, it follows from the Weierstrass approximation theorem that any such $Q$ can be uniformly approximated by sequences of polynomials on $[0,1]$ that satisfy (\ref{Q-funcopt}).



We now state our main result.

\begin{theorem}\label{precimain}
There exists $\theta_{0}>0$ such that whenever $\theta\in (0, \theta_{0})$, there exists $Q=Q_{\theta} \in C^{1}[0,1]$ satisfying (\ref{Q-funcopt}) such that $\kappa$ defined in (\ref{kappaprop}) satisfies $\kappa> 2\theta/3>0$.
\end{theorem}

The constant $\theta_{0}$ is explicitly computable. In fact, our numerical evidence (see Figure \ref{fig: diff_2third_act}) suggests the following explicit bound for any $0<\theta\le 1/2$:
\begin{align}\label{eqn: 2thirdbdd}
    \kappa \, > \, (2/3)\theta. 
\end{align}
Moreover, we will be able to make more general choices of $P$ and $Q$ in Section \ref{sect: gene}.


\subsection{Simplification}\label{sec: simplify}

Expanding the square of \eqref{eqn cpqroptim} and integrating the mixed term with the identity $\int_{0}^{1} \, f(x)f'(x) \, dx \, = \, (f(1)^2-f(0)^2)/2$, we have 
\begin{eqnarray} \label{eqn:cpr}
c(P,Q,R)=\frac{1}{2}+ J(Q),
\end{eqnarray}
where $J(Q)$ is the functional
\begin{eqnarray} \label{eqn:J} J(Q)=J_{R, \theta}(Q):= \int_0^1 \left(\frac{1}{\theta} w(y)^2 + \frac{\theta}{3}\,\, w'(y)^2\right)~dy,
\end{eqnarray}
and $w(y)$ is defined in \eqref{eqn cpqroptimwy}. Given constants $R, \theta>0$, we aim to find the minimum of $J(Q)$ subject to the condition \eqref{Q-funcopt}.

The expressions can be simplified by some changes of variables. Letting 
\begin{equation}\label{eqn:S}
   q(y):=Q(y+1/2) \qquad \text{and} \qquad S(t)=S_R(t):= q(t/2R)
\end{equation}
in \eqref{eqn:J}, we have
\begin{align}
    J(Q) \, = \,  \frac{e^R}{2R} \int_{-R}^R \bigg(\frac{1}{\theta} e^{t}S(t)^2 + \frac{\theta}{3}\,R^2 e^t \big( S(t)+2S'(t)\big)^2\bigg)~dt,
\end{align}
and the condition \eqref{Q-funcopt} becomes
\begin{align}\label{S-funcopt}
      S(-R)=1 \hspace{15pt} \text{and} \hspace{15pt} S(t) +S(-t) \, = \, 1.
   \end{align}
Here, $S$ is a real-valued, $C^{1}$-function defined on $[-R, R]$. It follows that $S(R)=0$, and, in particular, that
\begin{eqnarray*}\int_{-R}^R e^t S(t)S'(t)~dt \, = \,  -\frac{e^{-R}}{2}-\frac 12 \int_{-R}^R e^t S(t)^2~dt.
\end{eqnarray*}
Hence, we may write
\begin{eqnarray} \label{eqn:newJ} 
  J(Q)= -\frac{\theta}{3}\,R  \,  + \, \frac{e^R}{2R} \, K(S),
\end{eqnarray}
where 
\begin{eqnarray} \label{eqn:K0} K(S) = K_{\theta,R}(S):=\int_{-R}^R e^t\left[c_0S(t)^2 +c_1S'(t)^2)\right]~dt, \end{eqnarray}
and the constants $c_0,c_1$ are given by
\begin{equation}\label{eqn:c0c1}
 c_0 = c_0(R,\theta):= \frac{1}{\theta}-\frac{\theta}{3}\,R^2 \qquad \text{ and } \qquad  c_1= c_1(R,\theta):=\frac{4\theta}{3}\, R^2.
 \end{equation}
By splitting the integral in \eqref{eqn:K0} into the part from $-R$ to $0$ and the part from $0$ to $R$, we may impose the constraint $S(t)+S(-t)=1$ and further simplify $K(S)$ to write
\begin{eqnarray} \label{eqn:K}
    K(S)=  \int_{0}^R \Big(e^t \left( c_0S(t)^2+c_1S'(t)^2\right)+ e^{-t} \left( c_0(1-S(t))^2+c_1S'(t)^2\right)\Big) \, dt.
\end{eqnarray}
We have reduced the problem to finding the minimizer of $K(S)$ subject to \eqref{S-funcopt}. 

To keep our argument clean and reduce the number of variables in the course of proving Theorem \ref{precimain}, we introduce the following assumption:

\begin{assume}
Assume that $c_0  = c_1$, where $c_0,c_1$ are defined in \eqref{eqn:c0c1}.
 \end{assume}

We emphasize, however, that the steps leading up to Section \ref{sec: solveEL} do not depend on this assumption. A generalization is given in Section \ref{sect: gene}. The computations following Section \ref{sec: solveEL} can also be extended to the general case, though at the expense of increased computational complexity.

The assumption is equivalent to
 \begin{align}\label{simplipara}
     \theta=R^{-1}\sqrt{\frac{3}{5}}. 
 \end{align}
In this case, we have $c_0= c_1= 4/(5\theta)$. We note that $R$ can be taken to be any positive constant. We insert \eqref{simplipara} into \eqref{eqn:K} to write
\begin{eqnarray} \label{eqn:K1} 
K^{\ast}_{R}(S) \, := \,  \frac{4}{5\theta} \, \int_{0}^R e^t \left(S(t)^2+S'(t)^2 \right) +e^{-t} \left((1-S(t))^2+S'(t)^2 \right) ~dt.
\end{eqnarray}
Here, `$\ast$' in $K^{\ast}_{R}(S)$ signifies the use of (\ref{simplipara}). Our main problem can now be stated as follows. 
\begin{problem}\label{mainprob}
Given any $R>0$, find the minimum value of the functional $K_{R}^{\ast}(S)$ and the minimizer $S=S_{R}(t)$ subject to the condition \eqref{S-funcopt}.
\end{problem}
In the next section, we solve this problem by using the calculus of variations.


\section{Calculus of variations}\label{sect: COV}

In the following, we write $K^{\ast}(S)  = K^{\ast}_{R}(S)$. Suppose that $K^{\ast}(S)$ attains a minimum at $S$, and let $\epsilon>0$ be small and $\phi(t)$ be an arbitrary smooth test function on $[0,R]$ such that $\phi(0)=\phi(R)=0$. Then $K^{\ast}(S+\epsilon \phi)$, as a  function of $\epsilon$, has a minimum at $\epsilon =0$. Since 
\begin{eqnarray*}
    \frac{d}{d\epsilon}K^{\ast}(S+\epsilon \phi)\Big|_{\epsilon=0} &=&  \frac{4}{5\theta}\int_{0}^R \big(e^t  (2S \phi+2S'\phi') + e^{-t} ( 2 (S-1)\phi +2S'  \phi' )\big)~dt\\
     &=&
    \frac{16}{5\theta}\int_0^R\big( \cosh(t) S(t)-(S'(t) \cosh t)'-\frac{1}{2} e^{-t}\big) \phi(t) ~dt,
\end{eqnarray*}
we must determine $S$ for which the Euler--Lagrange equation
    \begin{eqnarray} \label{eqn:diffeq}
    \frac{d}{dt} \left(S'\cosh t\right)=S \cosh t -\frac{1}{2}e^{-t}
    \end{eqnarray}
 is satisfied. Equation (\ref{eqn:diffeq}) can be rewritten as
        \begin{eqnarray} \label{eqn:sdiffeq} S''(t) +(\tanh t) S'(t) - \big( S(t)-\frac{1}{1+e^{2t}}\big)=0.
  \end{eqnarray}
The boundary conditions for (\ref{eqn:sdiffeq}) are
   \begin{eqnarray} \label{eqn:init} 
   S(0)=1/2 \hspace{20pt} \text{and} \hspace{20pt} S(R)= 0. 
   \end{eqnarray} 

Before proceeding to solve for $S$, we incorporate the preceding steps to rewrite $c(P,Q,R)$ as defined in \eqref{eqn:cpr}. Assume that $S$ satisfies the differential equation (\ref{eqn:diffeq}). Then
   \begin{eqnarray*}
    \int_{0}^R S'(t)^2 \cosh t ~dt&=& S(t) S'(t)\cosh t\bigg|_{0}^R-  \int_{0}^R S(t)\frac d{dt}(S'(t)\cosh t)  ~dt\\
   &=& - \frac{1}{2} S'(0)-\int_{0}^R S(t)  \big( S(t) \cosh t - \frac{1}{2} e^{-t}\big) ~dt. 
   \end{eqnarray*}
It follows from \eqref{eqn:K} that
   \begin{align}\label{eqn:finalK}
   K^{\ast}(S) &= \frac{8}{5\theta} \int_{0}^R \big(  S(t)^2 \cosh t  + S'(t)^2 \cosh t -  e^{-t} S(t)\big) \, dt \, + \, \frac{4}{5\theta}(1-e^{-R})\notag \\
        &=\frac{4}{5\theta}\Big(1-e^{-R} - S'(0)  -  \int_{0}^Re^{-t} S(t)~dt \Big).
   \end{align} 
Upon substituting (\ref{simplipara}) and \eqref{eqn:finalK} into \eqref{eqn:newJ}, we find that \eqref{eqn:cpr} becomes
    \begin{align}\label{simplcPRQ}
 c(P,Q, R)&=\frac 1 2 -\frac{\theta}{3} R  +\frac{e^{R}}{2R} \frac{4}{5\theta} \Big(1-e^{-R} -  S'(0)  -   \int_{0}^Re^{-t} S(t)~dt\Big) \nonumber\\
&= \frac 1 2  + \frac{1}{\sqrt{15}} \, \Big(-1+2e^R\big(1-e^{-R} -  S'(0)  -  \int_{0}^Re^{-t} S(t)~dt\big)\Big).
\end{align} 
Readers should keep in mind that (\ref{simplcPRQ}) depends only on $R$ because of (\ref{simplipara}) and  (\ref{linearP}).


\subsection{Solving the Euler--Lagrange equation}\label{sec: solveEL}
 We first present the solution to the differential equation
\begin{align}\label{eqn: mainODEbdry}
    \begin{cases}
        S''(t)+\tanh(t)S'(t)-S(t)+(1+e^{2t})^{-1} =0, \\
        \hspace{40pt} S(0)=1/2, \hspace{10pt} \text{and} \hspace{10pt}  S(R)=0. 
    \end{cases}
\end{align}
It takes the form
\begin{align} \label{eqn: SRT}
   S_R(t):= C_1(R)f(t)+ g_0(t)+g_1(t) \int_0^t v_1(u)~du  +g_2(t) \int_0^t v_2(u)~du,
\end{align}
where the functions $C_1, f, g_0, g_1,g_2,v_1$ and $v_2$ are all expressible in terms of the Gauss ${}_2F_1$-hypergeometric function
\begin{align}\label{2F1ser}
      _{2}F_{1}(a,b, c;z) \, := \, \frac{\Gamma(c)}{\Gamma(a)\Gamma(b)} \sum_{n=0}^{\infty} \, \frac{\Gamma(a+n)\Gamma(b+n)}{\Gamma(c+n)} \frac{z^n}{n!}.
\end{align}
The series (\ref{2F1ser}) converges for $|z|<1$ and $a,b,c \in \mathbb{C}$ with $c\neq 0, -1, -2, \ldots$, and it admits an analytic continuation in the $z$-variable over $\mathbb{C}- [1,+\infty)$. We also record Euler's reflection formula, which states  
\begin{equation}\label{eqn:reflection}
\Gamma(s)\Gamma(1-s)=\frac{\pi}{\sin(\pi s)}, \qquad s \notin \mathbb{Z}.
\end{equation}

The component functions of $S_R(t)$ are defined as follows. To simplify the typesetting, we define $\phi:= (1+\sqrt{5})/2$;
\begin{align}
&\hspace{12pt} F^+(t):={}_2F_1\big(1/2, \, \phi, \, 1/2+\phi, \, -e^{2t}\big); \quad F^-(t):={}_2F_1\big(1/2, \, -\phi^{-1}, \, 1/2-\phi^{-1}; \, -e^{2t}\big), \nonumber\\
&F_1^+(u):={}_2F_1\big(1/2, \, 1+\phi, \, 1/2+\phi; \, -e^{2u}\big), \quad  F_1^-(u):={}_2F_1\big(3/2, \, -\phi^{-1}, \, 1/2-\phi^{-1}; \, -e^{2u}\big). \nonumber
\end{align}
The various components appearing on the right-hand side of \eqref{eqn: SRT} are given by:
\begin{align}
    g_1(t)= e^{-\phi^{-1}t}F^-(t), \hspace{5pt}
    g_2(t)= e^{\phi t}F^+(t), \hspace{5pt}  f(t) = g_1(t)- \frac{g_{1}(0)}{g_{2}(0)}\, g_2(t), \hspace{5pt} g_0(t)= \frac{1}{2}\frac{g_{2}(t)}{ g_{2}(0)}, \nonumber
\end{align}
\begin{align*}
v_1(u) \, = \, \frac{g_{2}(u)}{\mathcal{W}(u)(1+e^{2u})}, \hspace{15pt} 
v_2(u) \, = \,  - \ \frac{g_{1}(u)}{\mathcal{W}(u)(1+e^{2u})},
\end{align*}
and
\begin{align}\label{const: C1R}
    C_1(R)=\frac{-2 F^-(R)F^+(0)w_1(R)-e^{\sqrt{5}R}F^+(R)(1+2F^+(0)w_2(R))}
      {2F^-(R)F^+(0)-2 e^{\sqrt{5}R}F^-(0)F^+(R)},
\end{align}
where $w_1(t):=\int_0^t v_1(u) \, du$, $w_2(t):=\int_0^t v_2(u)\, du$, and
\begin{align}\label{spec: Wronsk}
    \mathcal{W}(u) \, := \,  e^{u}(2\phi F^-(u) F_1^+(u)-F^+(u) F_1^-(u)).
\end{align}
 Thus, $S_R(t)$ is fully described, and plots of its component functions are included. 

 We now derive \eqref{eqn: SRT}. Recall that the $_{2}F_{1}$-hypergeometric function can be characterized by a second-order linear homogeneous differential equation \cite[eq. (15.10.1)]{DLMF}. By a simple change of variables, one may readily observe (or check) that a set of fundamental solutions to the homogeneous differential equation $S''(t)+\tanh(t)S'(t)-S(t)=0$ is given by  $\{g_{1}(t), \ g_{2}(t)\}$.
Its Wronskian can be computed to be exactly $\mathcal{W}(u)$ defined in (\ref{spec: Wronsk}). 

The standard method of \textit{variation of parameters} (see \cite[Chapter 3]{CL55}) may now be used to solve our desired inhomogeneous equation (\ref{eqn: mainODEbdry}), which asserts that the sum of the last two terms on the right side of (\ref{eqn: SRT}) satisfies such a differential equation. The first two terms on the right side, which clearly satisfy the homogeneous equation, are included to ensure the boundary conditions are satisfied. Indeed, $f(0)=0$ and $g_{0}(0)=1/2$, and $C_{1}(R)$ is defined such that $S(R)=0$. Altogther, these yield a solution of the form (\ref{eqn: SRT}). Indeed, a direct \textit{Mathematica} check (\texttt{ODE\_Check.nb}) confirms that (\ref{eqn: SRT}) satisfies (\ref{eqn: mainODEbdry}).

 \begin{figure}
 \begin{subfigure}[b]{0.4\textwidth}
	\includegraphics[scale=0.75]{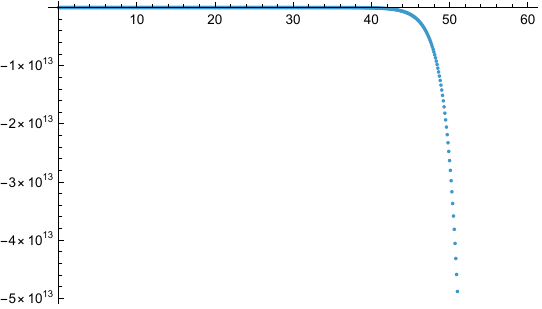}
	\caption*{Plot for $f(t)$}
    \end{subfigure}
    \hspace{3em}
    \begin{subfigure}[b]{0.4\textwidth}
	\includegraphics[scale=0.75]{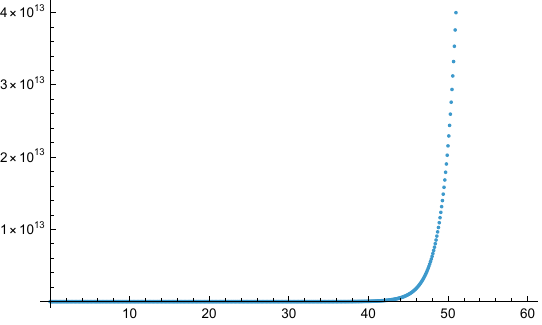}
	\caption*{Plot for $g_{0}(t)$}
	\label{fig: g0t}
    \end{subfigure}

    \vspace{2em}
    
    \begin{subfigure}[b]{0.4\textwidth}
	\includegraphics[scale=0.75]{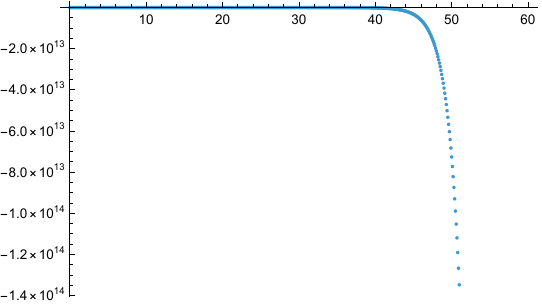}
	\caption*{Plot for $g_{1}(t)$}
	\label{fig: g1t}
    \end{subfigure}
        \hspace{3em}
    \begin{subfigure}[b]{0.4\textwidth}
	\includegraphics[scale=0.75]{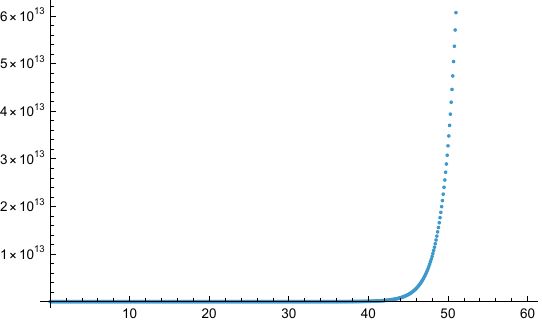}
	\caption*{Plot for $g_{2}(t)$}
	\label{fig: g2t}
    \end{subfigure}
\end{figure}

\begin{figure}
  \begin{subfigure}[b]{0.4\textwidth}
	\includegraphics[scale=0.75]{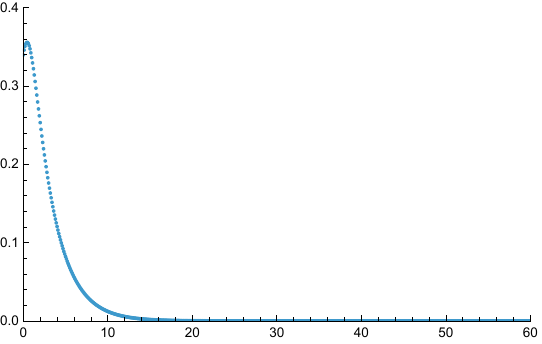}
	\caption*{Plot for $v_{1}(u)$ from 0 to 60}
	\label{fig: v1u}
    \end{subfigure}
    \hspace{3em}
    \begin{subfigure}[b]{0.4\textwidth}
	\includegraphics[scale=0.75]{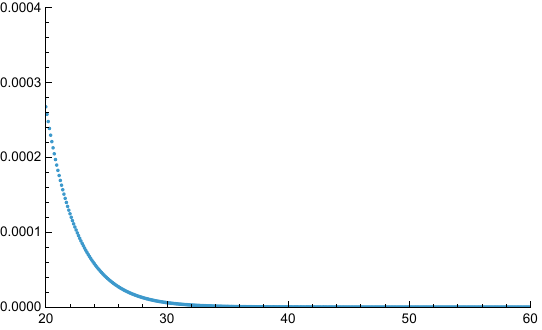}
	\caption*{Plot for $v_{1}(u)$ from 20 to 60}
	\label{fig: v1utail}
    \end{subfigure}
    
    \vspace{2em}
    
     \begin{subfigure}[b]{0.4\textwidth}
    \includegraphics[scale=0.75]{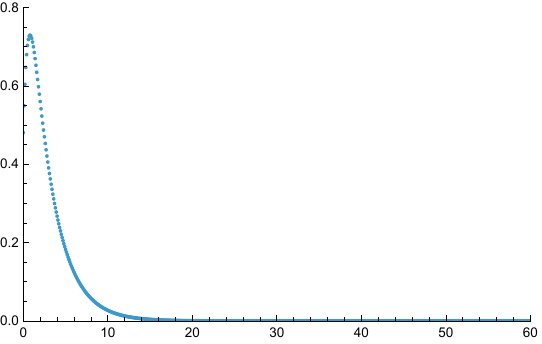}
	\caption*{Plot for $v_{2}(u)$ from 0 to 60}
	\label{fig: v2u}
    \end{subfigure}
    \hspace{3em}
    \begin{subfigure}[b]{0.4\textwidth}
    \includegraphics[scale=0.75]{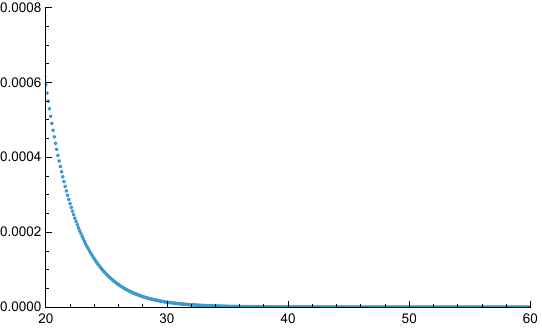}
	\caption*{Plot for $v_{2}(u)$ from 20 to 60}
	\label{fig: v2utail}
    \end{subfigure}

\end{figure}

\begin{figure}
\begin{subfigure}[b]{0.4\textwidth}
    \includegraphics[scale=0.75]{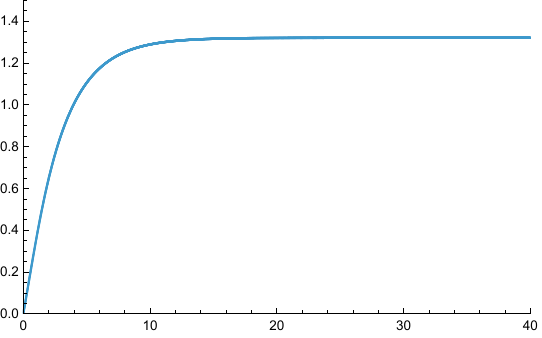}
	\caption*{Plot for $w_{1}(t)$ from 0 to 40}
	\label{fig: w1t}
    \end{subfigure}
    \hspace{3em}
    \begin{subfigure}[b]{0.4\textwidth}
    \includegraphics[scale=0.75]{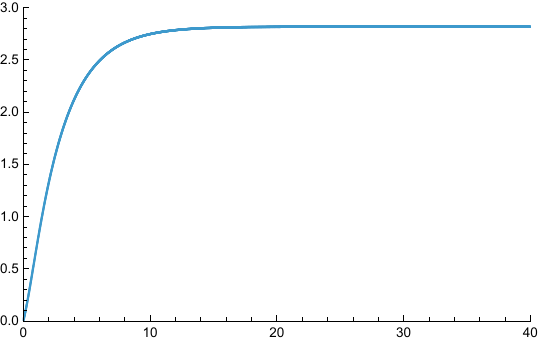}
	\caption*{Plot for $w_{2}(t)$ from 0 to 40}
	\label{fig: w2t}
    \end{subfigure}
\end{figure}

\begin{figure}
 \begin{subfigure}[b]{0.4\textwidth}
	\includegraphics[scale=0.75]{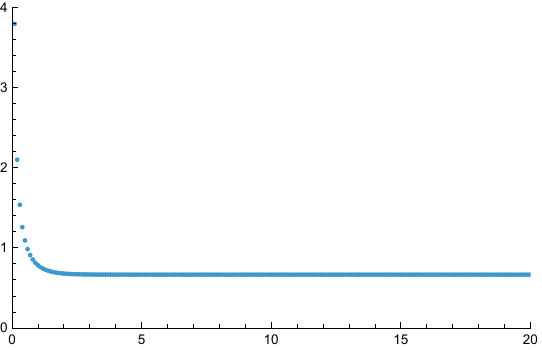}
	\caption*{Plot for $C_{1}(R)$}
    \end{subfigure}
    \hspace{3em}
    \begin{subfigure}[b]{0.4\textwidth}
    \includegraphics[scale=0.75]{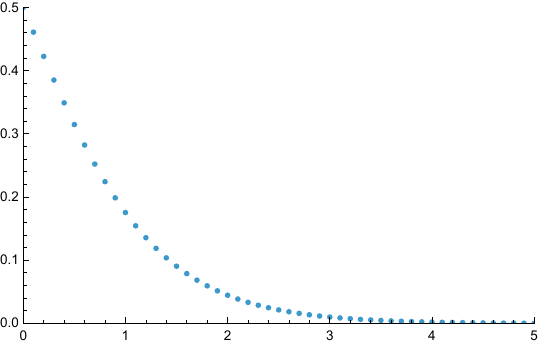}
	\caption*{Plot for $S_{R}(t)$ with $R=5$}
	\label{fig: S5t}
    \end{subfigure}

    \begin{subfigure}[b]{0.4\textwidth}
    \includegraphics[scale=0.75]{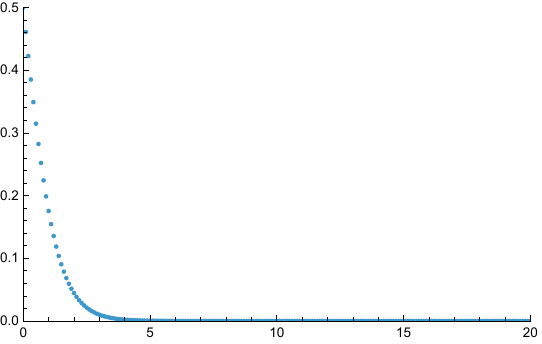}
	\caption*{Plot for $S_{R}(t)$ with $R=20$}
	\label{fig: S20t}
    \end{subfigure}
    \hspace{3em}
    \begin{subfigure}[b]{0.4\textwidth}
    \includegraphics[scale=0.75]{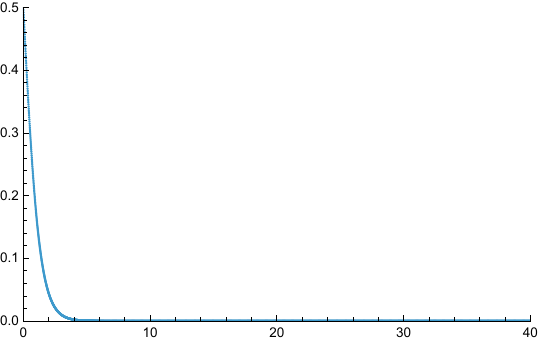}
	\caption*{Plot for $\lim_{R\to \infty} S_{R}(t)$}
	\label{fig: S10t}
    \end{subfigure}
\end{figure}


\subsection{Evaluating the components of $S_R(t)$ }\label{sec: asymps} We now record the leading order asymptotics of the above quantities  (as $t, u, R \to \infty$), computed with {\it Mathematica} (see the notebook \texttt{Asymp.nb} \footnote{ The commands used in the notebooks are also documented in Appendix \ref{sect: mathematica}.}). These asymptotics will be used in Section \ref{sec: proof} to complete the proof of Theorem \ref{precimain}. We first apply \cite[eq. (15.8.2)]{DLMF}, which for real $z<0$ states
\begin{align}\label{DLMFeq}
\frac{\sin\pi(b-a)}{\pi\Gamma(c)}  \,  & _2F_1(a,b,c;z) \nonumber\\
= \, & \frac{(-z)^{-a}}{\Gamma(b)\Gamma(c-a)\Gamma(a-b+1)} \, _2F_1\left(a, \, a-c+1, \, a-b+1; \, \frac 1 z\right) \nonumber\\
 &\hspace{5pt} -\,\frac{(-z)^{-b}}{\Gamma(a)\Gamma(c-b)\Gamma(b-a+1)}
    \, _2F_1\left(b, \, b-c+1, \, b-a+1;\, \frac 1 z\right).
\end{align}
In our setting, $z=-e^{2t}$, where $t\ge 0$, and $a,b$, and $c$ are non-integer real numbers. Thus by \eqref{DLMFeq}, the power series (\ref{2F1ser}), and the reflection formula \eqref{eqn:reflection}, we have the full expansion
\begin{align}\label{fullasymp}
    \frac{\sin \pi(b-a)} {\Gamma(c)}  \, & _2F_1(a,b,c; -e^{2t}) \nonumber\\
    = \, &  \frac{  e^{-2ta} \sin \pi (c-a)}{\Gamma(a)\Gamma(b)   }   \sum_{n=0}^{\infty} \, \frac{\Gamma(a+n)\Gamma(a-c+1+n)}{\Gamma(a-b+1+n)} \frac{(-e^{-2t})^n}{n!}\nonumber\\
    &\hspace{5pt} -
    \frac{e^{-2tb} \sin \pi(c-b)}{\Gamma(a) \Gamma(b)} \sum_{n=0}^{\infty} \, \frac{\Gamma(b+n)\Gamma(b-c+1+n)}{\Gamma(b-a+1+n)} \frac{(-e^{-2t})^n}{n!}. 
\end{align}
Note that \eqref{fullasymp} implies, as $t \to \infty$, that ${}_2F_1(a,b,c; -e^{2t}) = O(e^{-2t\min\{a,b\}})$. Using \eqref{fullasymp}, we may use {\it Mathematica} to determine the leading-order asymptotics for the functions making up $S_R(t)$  efficiently and precisely. First, as $t\to \infty$, we have \footnote{ In this article, the decimal values given are truncations of the actual values.  }
\[
F^+(t)\sim 1.2427\ldots e^{-t}\qquad \mbox{and} \qquad F^-(t)\sim-2.75957 \ldots e^{(\sqrt{5}-1)t}
\]
so that
\[
g_1(t)\sim-2.75957 \dots e^{\frac{\sqrt{5}-1}{2}t}, \qquad g_2(t)\sim 1.2427\dots e^{\frac{\sqrt{5}-1}{2}t},
\]
\vspace{5pt}
\[
\hspace{5pt}f(t)\sim - e^{\frac{\sqrt{5}-1}{2}t} \hspace{40pt} \text{ and } \qquad  g_0(t)\sim 0.81855\dots  e^{\frac{\sqrt{5}-1}{2}t}. 
\]

\noindent Next,  as $u\to \infty$, we have
\[
F_1^+(u)\sim 0.85875 \dots e^{-u} \qquad \text{and} \qquad F_1^-(u)\sim -6.17059\dots e^{(\sqrt{5}-1)u}
\]
so that
\[
v_1(u)\sim 0.55579\dots e^{-\frac{3-\sqrt{5}}{2}u} \qquad \text{ and }  \qquad v_2(u)\sim  \ 1.23411 \dots e^{-\frac{3-\sqrt{5}}{2}u}.
\]
It is now clear that $\lim_{R\to \infty} w_{i}(R)$ exist for $i=1,2$. 

\begin{remark}
    The coefficients appearing in the asymptotics above can indeed be expressed in closed form. For example, we have $1.2427\ldots = \Gamma(1+\sqrt{5}/2)\Gamma(\sqrt{5}/2)/\Gamma((1+\sqrt{5})/2)^2$ and $-2.75957 \ldots= \csc(\sqrt{5}\pi/2)$, which follow readily from the Pfaff transformation formula \cite[eq. (15.8.1)]{DLMF}, the Gauss summation formula \cite[eq. (15.4.20)]{DLMF}, and the Euler reflection formula (\ref{eqn:reflection}). 
\end{remark}

\vspace{5pt}

\noindent Finally, we have that the numerator and denominator of $C_{1}(R)$ are both $\asymp e^{(\sqrt{5}-1)R}$, and thus the limit of  $C_1(R)$ as $R\to \infty$ exists. Indeed,  by \textit{Mathematica}, we find that 
\begin{align}
    w_{1}(\infty)\,  = \, 1.3208 \ldots \hspace{20pt} \text{and} \hspace{20pt} w_{2}(\infty) \, = \, 2.8166\ldots.
\end{align}
Thus,
\begin{equation}\label{eqn: C1limit}
C_1(R) \, \, \sim \, \frac{4.18979\ldots w_{1}(\infty)-1.2428\ldots-1.88691\ldots w_{2}(\infty)}{-1.51827\ldots} \, \sim \,  0.674\dots. 
\end{equation}

\noindent Once again, the constants $4.18979\ldots$, $1.2428\ldots$, $1.88691\ldots$, and $1.51827\ldots$ appearing in (\ref{eqn: C1limit}) can be expressed in terms of the $\Gamma$-functions and the $_{2}F_{1}$-hypergeometric functions evaluated at $z=-1$. Also, we note that the function $\lim_{R\to \infty} S_{R}(t)$ is of the order $e^{\frac{\sqrt{5}-1}{2}t}$. 

\begin{remark}
For the readers' convenience, we have generated tables of values for: 
\begin{itemize}
    \item $v_{1}(u)$, $v_{2}(u)$ over the interval $u\in [0,40]$ with  step size  $=1/5000$ \\
   \textnormal{(files: \texttt{Tbvv140oneOver5000.mx},  \texttt{Tbvv240oneOver5000.mx});} \vspace{4pt}

    \item $w_{1}(t)$, $w_{2}(t)$ over the interval $t\in [0,40]$ with  step size $=1/100$\\
    \textnormal{(files: \texttt{Tbww1t40ref1over100.mx},  \texttt{Tbww2t40ref1Over100.mx});} \vspace{4pt}

    \item $C_{1}(R)$ over the interval $R\in [0,40]$ with step size $=1/100$ \\
    \textnormal{(file: \texttt{ConstC40oneOver100.mx});} \vspace{4pt}

    \item $\lim_{R\to \infty} S_{R}(t)$ over the interval $t\in [0,40]$ with step size  $=1/100$ \\
    \textnormal{(file: \texttt{Sinf40oneOver100.mx})}. 
\end{itemize}
These datasets can be generated by the notebook \textnormal{\texttt{Pre\_compute.nb}}. The graphs included can be created with the standard plotting functions \textnormal{\texttt{Plot[]}} and \textnormal{\texttt{ListPlot[]}} (or see \textnormal{\texttt{Plot.nb}}). 
\end{remark}



\section{Relation to Siegel's function: Proof of Proposition \ref{appSige}}

  Recall from (\ref{eqn:S}) that $Q_{R}(y)=S_{R}(2R(y-1/2))$. The initial condition in (\ref{eqn: mainODEbdry}) implies that $Q_{R}(1/2)=0$. By (\ref{S-funcopt}), it suffices to consider the range $1/2<y\le 1$. 
  
  Let $y_{0}\in (1/2,1]$ be given. We have determined that $f(t)$, $g_{0}(t)$, $g_{1}(t)$, $g_{2}(t)$ all exhibit the same asymptotic order  $e^{\frac{\sqrt{5}-1}{2}t}$ as $t\to\infty$. Then for any $\mathcal{F} \in \{f, g_{0}, g_{1}, g_{2}\}$, we have
  \begin{align}
      \lim_{R\to\infty} \, \mathcal{F}(t)e^{-\frac{\sqrt{5}-1}{2}t} \Big|_{t=2R(y_{0}-1/2)} = \lim_{R\to\infty} \, \mathcal{F}(t)e^{-\frac{\sqrt{5}-1}{2}t} \Big|_{t=R}.  
  \end{align}
  Also, it is clear that $\lim_{R\to\infty} \, w_{i}(2R(y_{0}-1/2)) = \lim_{R\to \infty} w_{i}(R) = w_{i}(\infty)$ for $i=1,2$. Now, using (\ref{eqn: SRT}), (\ref{eqn: C1limit}) and the initial condition, we conclude that
   \begin{align}
       \lim_{R\to \infty} \, S_{R}(t)e^{-\frac{\sqrt{5}-1}{2}t}  \Big|_{t=2R(y_{0}-1/2)} \, = \,  \lim_{R\to \infty} \, S_{R}(R) e^{-\frac{\sqrt{5}-1}{2}R} \, = \, 0,
   \end{align}
   and thus, $\lim_{R\to \infty} \, Q_{R}(y_{0})=0$. This completes the proof.

\begin{figure}[H]
 \begin{subfigure}[b]{0.45\textwidth}
\includegraphics[scale=0.8]{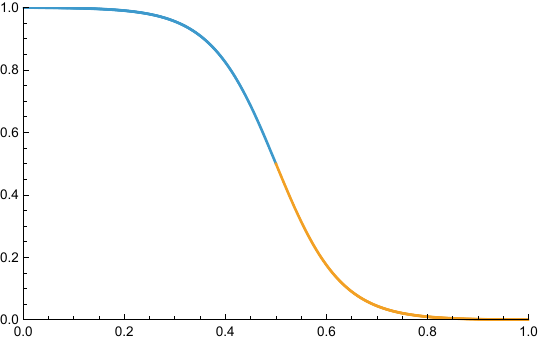}
\caption{Plot for $Q_{5}(y)$}
	\label{fig: Q5y}
\end{subfigure}
    \hspace{1.5em}
    \begin{subfigure}[b]{0.45\textwidth}
\includegraphics[scale=0.8]{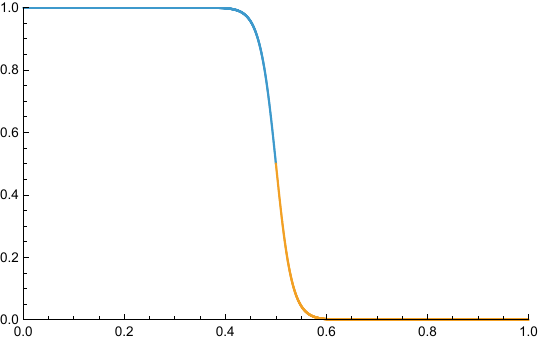}
\caption{Plot for $Q_{20}(y)$}
	\label{fig: Q20y}
\end{subfigure}
    \end{figure}


\section{Proof of Theorem \ref{precimain}}\label{sec: proof}

Our goal is to show that
\begin{equation}\label{eqn: kappaupdated}
\kappa \, = \, \kappa(\theta) \, := \,  1 - \frac{\log c(P,Q,R)}{R} >0,
\end{equation}
where $\theta=R^{-1}\sqrt{3/5},  $ $P(x)=x$, $Q(y)=Q_{R}(y)=S_{R}(2R(y-1/2))$, and 
\begin{equation}\label{eqn: crRdef}
 c(P,Q,R) \, = \, \frac 1 2 - \frac{1}{\sqrt{15}} +\frac{2}{\sqrt{15}}e^R\Big(1-e^{-R} -  S_R'(0)  -  \int_{0}^Re^{-t} S_R(t)~dt\Big),
\end{equation}
which depends only on $R$ (or equivalently, $\theta$). In particular, we seek an upper bound for the expression enclosed in $(\cdots)$ in (\ref{eqn: crRdef}).

The asymptotics given in Section \ref{sec: asymps} imply that the following integrals converge absolutely. Moreover, their numerical values can be computed with \texttt{Integral\_kappa.nb}: 
\[\int_0^\infty e^{-t} f(t) ~dt = -2.1166\dots, \quad \quad \int_0^\infty e^{-t} g_0(t) ~dt =  1.9453\dots,\]
\vspace{5pt}
\[\int_0^\infty e^{-t} g_1(t) \int_0^t v_1(u)~du~dt=-4.294\dots, \quad \quad \int_0^\infty e^{-t} g_2(t) \int_0^t v_2(u)~du~dt= 4.024\dots.\]
From these and (\ref{eqn: SRT}), it follows that
\begin{align}\label{eqn: finSRT}
    \int_0^R e^{-t}S_R(t)~dt \, = \,  -2.1166\dots C_{1}(R) \, + \, 1.675\ldots, 
\end{align}
and
\begin{equation}\label{eqn: limSRT}
\lim_{R\to \infty}\int_0^R e^{-t}S_R(t)~dt \ge 0.248.
\end{equation}

Also, we have the numerical values
\begin{align}
    &f'(0)=-1.47277\dots, \hspace{15pt} g_0'(0)=0.602775\dots, \hspace{15pt} g_1(0)=-1.07479\dots, \nonumber\\
    &g_2(0)=0.759136\dots, \hspace{15pt} v_1(0)=0.339496\dots, \hspace{15pt} v_2(0)=0.480664\dots.
\end{align}
By \eqref{eqn: SRT}, (\ref{eqn: C1limit}) and the evaluation above,  we obtain
\begin{align}\label{eqn: finderSP}
    S_R'(0)&= C_1(R) f'(0)+ g_0'(0)+g_1(0)  v_1(0) +g_2(0)   v_2(0) \nonumber\\
    \, &= \, -1.47277\ldots C_{1}(R) \, + \, 0.60277\ldots
\end{align}
and
\begin{equation}\label{eqn: limSP}
\lim_{R\to \infty} \, S_R'(0) \, \ge \,  -0.39006. 
\end{equation}

Therefore, from (\ref{eqn: finSRT}) and (\ref{eqn: finderSP}), we have
\begin{align}
    1-e^{-R} -  S_R'(0)  -  \int_{0}^Re^{-t} S_R(t)~dt \, = \,  -e^{-R}+3.58938\ldots C_{1}(R) \, - \, 1.278\ldots.
\end{align}
By \eqref{eqn: limSRT} and \eqref{eqn: limSP}, we have
\[
\lim_{R\to \infty}\big(1-e^{-R} -  S_R'(0)  -  \int_{0}^Re^{-t} S_R(t)~dt\big) \, \le \, 1.1421. 
\]

Since $R=\theta^{-1}\sqrt{3/5}$, there exists $\theta_0>0$ such that for all $0 < \theta < \theta_0$, we have
\[
1-e^{-R} -  S_R'(0)  -  \int_{0}^Re^{-t} S_R(t)~dt \, \le \,  1.2, 
\]
and
\begin{align*}
 c(P, Q, R) \, \le \,  \frac{1}{2}  +\frac{1}{\sqrt{15}}(2.2842 e^{R}-1) \, \le \, 0.5898e^{R} + 0.242 \, \le \, 0.5898e^{R}(1+0.411e^{-R}).
\end{align*}
Using the elementary inequality $\log (1+x)\le x$ for $x>-1$, it follows that
\begin{align}
    \kappa \, \ge \, \frac{1}{R}\big(-\log 0.5898 - 0.411 e^{-R}\big).
    \end{align}
There exists $\theta_{1}>0$ such that whenever $0<\theta<\theta_{1}$ we have $0.411 e^{-R}< 0.0001$. Thus,  
\begin{align}
     \kappa \, \ge \, \theta \sqrt{\frac{5}{3}}(-\log 0.5898 -0.0001) \, \ge \, 2\theta/3. 
\end{align}
for $0<\theta<\min\{\theta_{0}, \theta_{1}\}$. In particular, there is a positive proportion of zeros for $\zeta(s)$ for sufficiently small $\theta>0$.  This completes the proof of Theorem \ref{precimain}.



     \section{Some data}
We record the proportions $\kappa=\kappa(\theta)$ of zeros of $\zeta(s)$ on $\re s=1/2$ using (\ref{kappaprop}), (\ref{simplcPRQ}), and (\ref{eqn: SRT}) for various values of $\theta$. Here are the plots and numerics (generated with \texttt{Integral\_kappa.nb}):
\begin{figure}[H]
 \begin{subfigure}[b]{0.45\textwidth}
\includegraphics[scale=0.8]{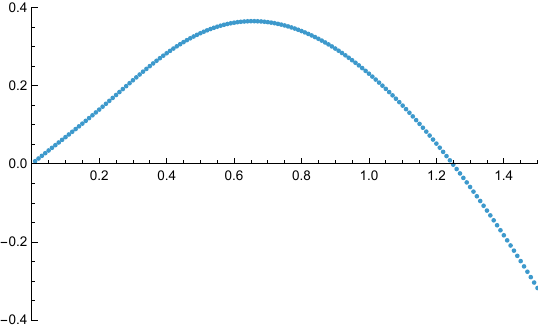}
\caption{Plot for $\kappa(\theta)$ with $P(x)=x$, $Q(x)=Q_{R}(x)$, $R=\theta^{-1}\sqrt{3/5}$}\label{fig: kappaPlot_Lev_Op}
\end{subfigure}
    \hspace{1.5em}
    \begin{subfigure}[b]{0.45\textwidth}
\includegraphics[scale=0.9]{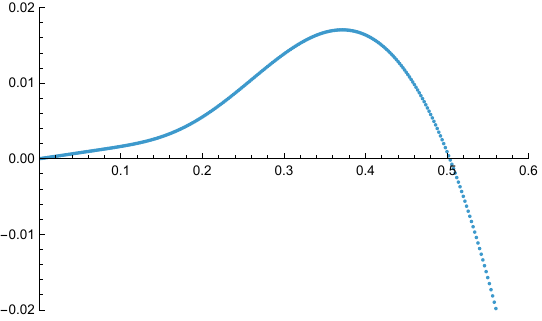}
	\caption{Plot for $\kappa(\theta)-(2/3)\theta$ }
	\label{fig: diff_2third_act}
    \end{subfigure}   
    \end{figure}
\vspace{-2em}
    
  \begin{eqnarray*}
     \begin{array}{|c|c|}
     \hline
     \theta& \kappa\\
     \hline
     2/3 & 0.364\dots\\
     \hline
     1/2 & 0.334\dots\\
     \hline
     1/4 & 0.176\dots\\
     \hline
     1/6 &  0.114\dots\\
     \hline
     1/8& 0.0854\dots\\
     \hline
     5/54 = 0.092\ldots& 0.0632\dots\\
     \hline
     1/100 & 0.00682\dots\\
     \hline
     1/500 & 0.00136\dots\\
     \hline
    \end{array}
    \end{eqnarray*}

The values  $5/54$ and $1/8$ in our table for the $\zeta$-function are of interest as they correspond to $5/27$ and $1/4$ for modular $L$-functions in \cite{Ber15} and \cite{KRZ19}, based on the well-known principle that a mollifier of length $\theta$ for a $\mathrm{GL}(n)$ $L$-function corresponds to a mollifier of length $\theta/n$ for a $\mathrm{GL}(1)$ $L$-function; see \cite[pp. 215-216]{Far94} and \cite[Theorem 5]{Ber15}.   The proportions of 0.0297 and 0.0297607 were obtained unconditionally in \cite{Ber15} and \cite{KRZ19}, respectively, and 0.0693 and 0.0693872  assuming the Ramanujan--Selberg conjecture (RSC). (It is worth noting that a two-piece mollifier was used in \cite{KRZ19}.) In our table, we obtain 0.0632 unconditionally and 0.0854 on RSC. Using stronger spectral inputs on shifted convolution sums and a more extensive computer search, the best-known proportions for modular $L$-functions are 0.0696 unconditionally and 0.0896 on RSC (\cite{AT21}). 


\section{Generalizations and further questions}\label{sect: gene}


We sketch the general set-up.  Let  $P\in C^{1}[0,1]$ be a \textit{given} real-valued function such that $P(0)=0$ and $P(1)=1$,
and we let 
   $$B:=\int_0^1P(x)^2 ~dx \qquad \mbox{and} \qquad C:=\int_0^1P'(x)^2~dx.$$ 

Let $\beta\in \mathbb{R}$ and $R>0$ be constants at our disposal. By the same reasoning as in Section \ref{sect: COV}, we are readily led to the problem of optimizing the functional
\begin{eqnarray}\label{eqn:Kfuncgen}
    K(S) \, := \,  \int_{0}^R \Big(e^t ( c_0S(t)^2+c_1S'(t)^2)+ e^{-t} ( c_0(\beta-S(t))^2+c_1S'(t)^2)\Big) \, dt,
    \end{eqnarray}
with $S\in C^{1}[0, R]$ and subject to boundary conditions 
\begin{eqnarray} \label{eqn: geninit} 
   S(0) \, = \, \beta/2 \hspace{20pt} \text{and} \hspace{20pt}   S(R) \, = \, \beta-1, 
   \end{eqnarray} 
 where
    \begin{align*}
       c_0\, := \,  \frac{C}{\theta}-\theta BR^2 \hspace{20pt} \text{and} \hspace{20pt} c_1 \, := \, 4\theta B R^2. 
   \end{align*}
For any $S$ satisfying the constraints above and for any admissible value of $\theta>0$ (cf. (\ref{Conmollif})), we have
   \begin{align}\label{kappa: general}
       \kappa \, = \, 1-\frac1R \, \log \Big(\frac{1+e^{2R}(\beta-1)^2}{2} \, + \, c_{1}\, \frac{e^{R}}{2R}\frac{e^R (1-\beta)^2-e^{-R}}{2} \, + \, \frac{e^{R}}{2R}\,  K(S)\Big). 
   \end{align}
 Furthermore, for the optimal function $S$, we may simplify (\ref{eqn:Kfuncgen}) as
   \begin{eqnarray} \label{eqn: generfinalK}
   K(S)=c_0\beta^2(1-e^{-R}) -c_1 \beta S'(0) +2(\beta-1)c_1S'(R)\cosh R -c_0 \beta  \int_{0}^Re^{-t} S(t)~dt.
   \end{eqnarray}  

\begin{remark}
The constant $\beta$ arises from the constraint for (real-valued) $Q\in C^{1}[0,1]$:  \begin{align}
        Q(0)=1 \hspace{15pt} \text{and} \hspace{15pt} Q(y) +Q(1-y) \, = \, \beta. 
   \end{align}
   Again, we have $S(-t):= \beta-S(t)$ for $t\ge 0$ and $Q(y)=S(2R(y-1/2))$ for $y\in [0,1]$.
\end{remark}

We also let 
\begin{eqnarray*}
   c \, := \, -\frac {c_0}{c_1} \, = \, \frac 1 4 -\frac{C}{4B\theta^2 R^2}. 
   \end{eqnarray*}
With (\ref{eqn: geninit}), the Euler--Lagrange equation for $K(S)$ in this general case is given by
\begin{eqnarray}
        S''(t) \, + \, (\tanh t) S'(t) \, + \, c  S(t) \, = \, \frac{c\beta}{1+e^{2t}}.
  \end{eqnarray}

\begin{remark}
    The Euler--Lagrange equation above depends on $R>0$, $\beta\in\mathbb{R}$, and $c<1/4$ which further depends on $B, C, \theta, R$. When proving Theorem \ref{precimain}, the only parameter that varies is $R$ (or equivalently, $\theta$). In that case, $c=-1$, $\beta=1$, $B=1/3$, and $C=1$. 
\end{remark}

Let $\phi_{c}:=(1 + \sqrt{1-4 c})/2$. A set of fundamental solutions to the homogeneous equation of (\ref{eqn:sdiffeq}) is given by 
\begin{align}
    g_{1}(t; c) \, &:= \, \exp(tc/\phi_{c}) \, _{2}F_{1}\big(1/2, \, c /\phi_{c}; \, 1/2- c/\phi_{c}; \, -e^{2t} \big), \nonumber\\
    g_{2}(t; c) \, &:= \, \exp(t\phi_{c}) \, _{2}F_{1}\big(1/2, \, \phi_{c}; \, 1/2 +\phi_{c}; \, -e^{2t} \big). 
\end{align}   
Parallel to the notations in Section \ref{sec: solveEL}, we set
\begin{align}
&\hspace{12pt} F^+(t;c):={}_2F_1\big(1/2, \, \phi_{c}, \, 1/2+\phi_{c}; \, -e^{2t}\big); \quad F^-(t;c):={}_2F_1\big(1/2, \, c/\phi_{c}, \, 1/2+c/\phi_{c}; \, -e^{2t}\big); \nonumber\\
&F_1^+(u;c):={}_2F_1\big(1/2, \, 1+\phi_{c}, \, 1/2+\phi_{c}; \, -e^{2u}\big); \quad  F_1^-(u;c):={}_2F_1\big(3/2, \, c/\phi_{c}, \, 1/2+c/\phi_{c}; \, -e^{2u}\big). \nonumber
\end{align}
The Wronskian for $\{g_{1}(t; c), \, g_{2}(t; c)\}$ can be computed to be
\begin{align}
    \mathcal{W}(u;c) \, = \, \mathcal{W}(g_{1},g_{2})(u;c) \, = \, e^{u}(2\phi_{c} F^-(u;c) F_1^+(u;c)-F^+(u;c) F_1^-(u;c)). 
\end{align}
Similar to before, we introduce:
\begin{align}
    f(t;c) = g_1(t; c)- \frac{g_{1}(0;c)}{g_{2}(0;c)}\, g_2(t;c), \hspace{15pt} g_0(t;c)= \frac{1}{2}\frac{g_{2}(t;c)}{ g_{2}(0;c)}, 
\end{align}
\begin{align}
    v_{1}(u; c) \, := \,  \frac{g_{2}(u;c)}{ \mathcal{W}(u;c)} \frac{1}{1+e^{2u}}, \hspace{15pt}  v_{2}(u; c) \, := \,  -\, \frac{g_{1}(u;c)}{ \mathcal{W}(u;c)} \frac{1}{1+e^{2u}}.
\end{align}
It follows that the solution to the Euler--Lagrange equation is given by
\begin{align}\label{eqn: genoptSRt}
    S_{R}(t; c, \beta) \, = \, C_1(R; c, \beta) f(t;c)+ \beta g_0(t;c) \, - \, c\beta  g_1(t;c) \int_0^t v_1(u;c)~du  - c\beta g_2(t) \int_0^t v_2(u;c)~du,
\end{align}
where 
\begin{align}
    C_{1}(R; c,\beta) \, = \, \frac{\beta-1- \beta g_0(R;c)+ c\beta  g_1(R;c)w_{1}(R;c)+c\beta g_{2}(R; c)w_{2}(R;c)}{f(R; c)}. 
\end{align}

We conclude this article with a few questions intended to stimulate future research.

\begin{enumerate}
\item In \cite{Con89}, the following choice of $P$ in the mollifier was used:
\begin{align}\label{optimalP}
    P(x) \, = \,  P_r(x) \, := \, \frac{\sinh (rx)}{\sinh(r)},
\end{align}
which was shown to be optimal. In this case, we have
\begin{align*}
   \hspace{25pt}  B=\frac{\coth (r)-r \text{csch}^2(r)}{2 r} \hspace{15pt} \text{and} \hspace{15pt} C=\frac{r}{2} \, \text{csch}(r) (\cosh (r)+r \text{csch}(r)).
\end{align*}
It is natural to ask what the largest possible value of $\kappa= \kappa(\theta)$, as given in (\ref{kappa: general}), would be when combining \eqref{optimalP} with our constructions \eqref{eqn: generfinalK} and \eqref{eqn: genoptSRt}. More precisely, for which choices of $\beta$, $r$, $R$ optimize the value of $\kappa$? \vspace{5pt}

\item Study the optimization problem (\ref{eqn:Kfuncgen}) \textit{jointly} in the variables $P$ and $S$. It is likely that $P$  and $S$ given in (\ref{optimalP}) and (\ref{eqn: genoptSRt}) satisfy the corresponding Euler--Lagrange equation. \vspace{5pt}

\item Our optimal linear combination, together with the optimal (Levinson-type) mollifier, should give better values of $\kappa$. When $\theta$ is large, say $\theta = 4/7-\epsilon$, it is of interest to know how much improvement this would offer over \cite{Con89}. 

We expect that the gain may be modest, likely not matching the state-of-the-art improvements achieved through more sophisticated mollifiers with multiple pieces, e.g., those developed in \cite{BCY11}, \cite{Fen12} and \cite{PR+20}. \vspace{5pt}


\item An interesting question is whether one can choose a simpler function $Q$ to prove our Theorem \ref{precimain}, albeit at the expense of weaker numerical constants for $\kappa$. \vspace{5pt}

\item Prove that $\kappa(\theta)> 2\theta/3$ \textit{for any} $\theta \in (0,1/2]$, using $P(x)=x$ and $Q(y)=Q_{R}(y)=S_{R}(2R(y-1/2))$ with $S_{R}$ defined in (\ref{eqn: SRT}).  

It suffices to establish a good lower bound for $C_{1}(R)$ for $R\ge 2\sqrt{3/5}$. \vspace{5pt}

\item In view of Proposition \ref{appSige},  we wonder, how does the configuration of zeros for the function $Q_{R}(-\frac{1}{L} \frac{d}{ds})\zeta(s)$ changes as $R$ varies, where $Q_{R}(y)=S_{R}(2R(y-1/2))$ is the `optimal' function determined in Section \ref{sec: solveEL}.  

While $Q_{R}(-\frac{1}{L} \frac{d}{ds})\zeta(s)$ serves as (pointwise) approximations to Siegel's function $\mathfrak{f}(s)$, we suspect that, as $R\to \infty$, its zeros are pushed to the left relatively far from the critical line.  This should stand in contrast to the distribution of zeros of $\mathfrak{f}(s)$ (see \cite[Fig. 5]{Con16}), which appears to exhibit clusters of zeros near the critical line, both to its left and right.

Below is a contour plot of $Q(-\frac{1}{L} \frac{d}{ds})\zeta(s)$, where $Q(y)=0.492+0.602(1-2x)-0.08(1-2x)^3-0.06(1-2x)^5+0.046(1-2x)^7$ is the polynomial used in \cite{Con89}. 
\end{enumerate}

\begin{figure}[H]
	\includegraphics[scale=0.95]{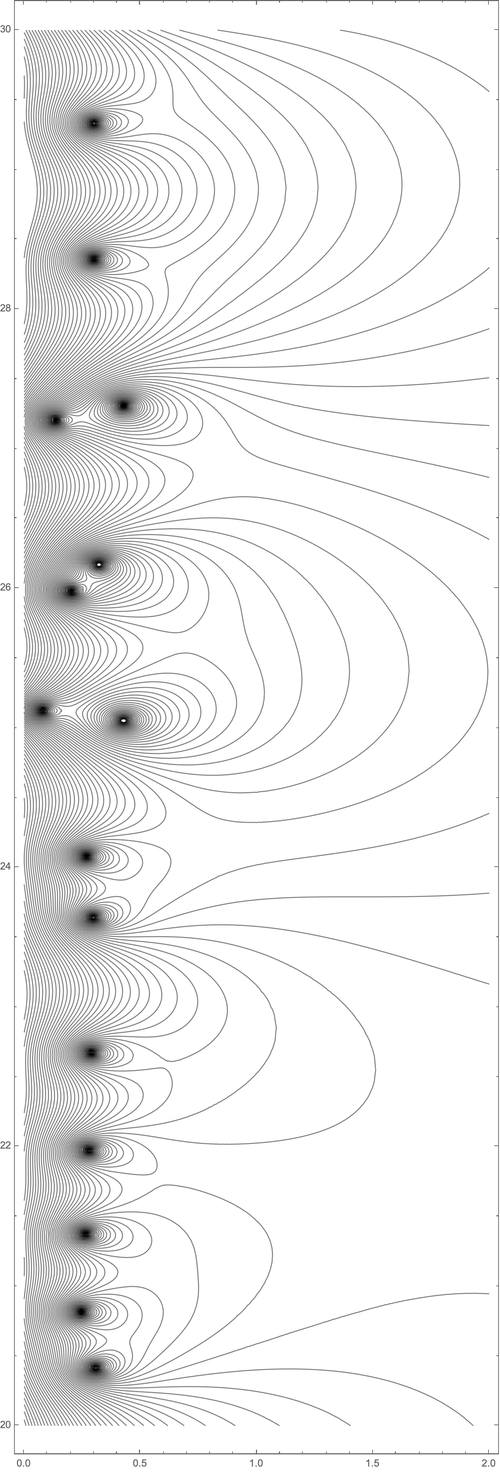}
	\label{fig: Qplot}
\end{figure}

\appendix

\section{Mathematica commands (\texttt{ODE\_Check.nb}, \texttt{Pre\_compute.nb}, \texttt{Asymp.nb}, \texttt{Integral\_kappa.nb})}\label{sect: mathematica}

\begin{center}
    \includegraphics[page=1, width=0.85\textwidth, trim=1.8cm 5cm 2cm 2.5cm, clip]{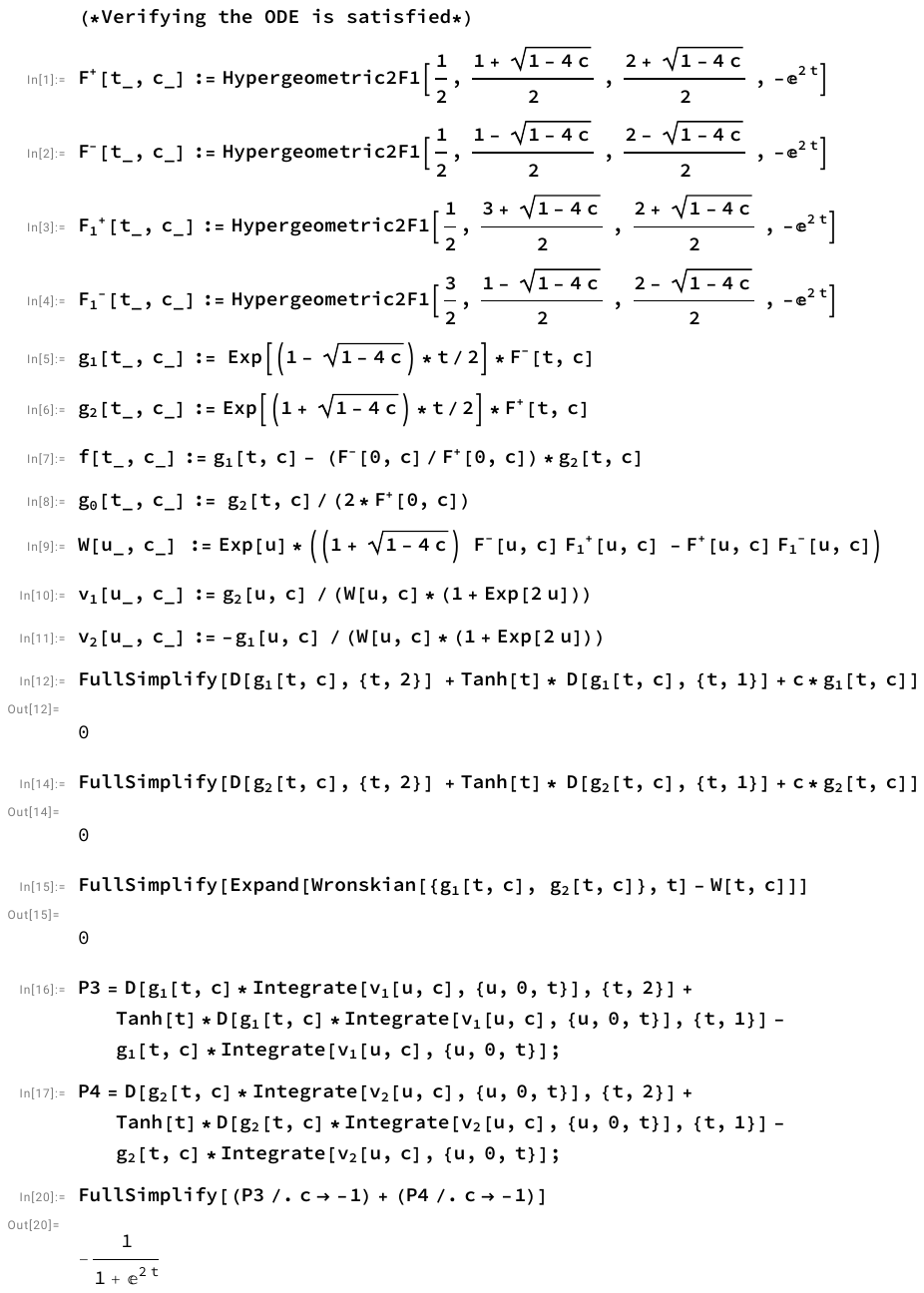}
\end{center}

\begin{center}
    \includegraphics[page=1, width=0.9\textwidth, trim=1.8cm 2cm 2cm 2.5cm, clip]{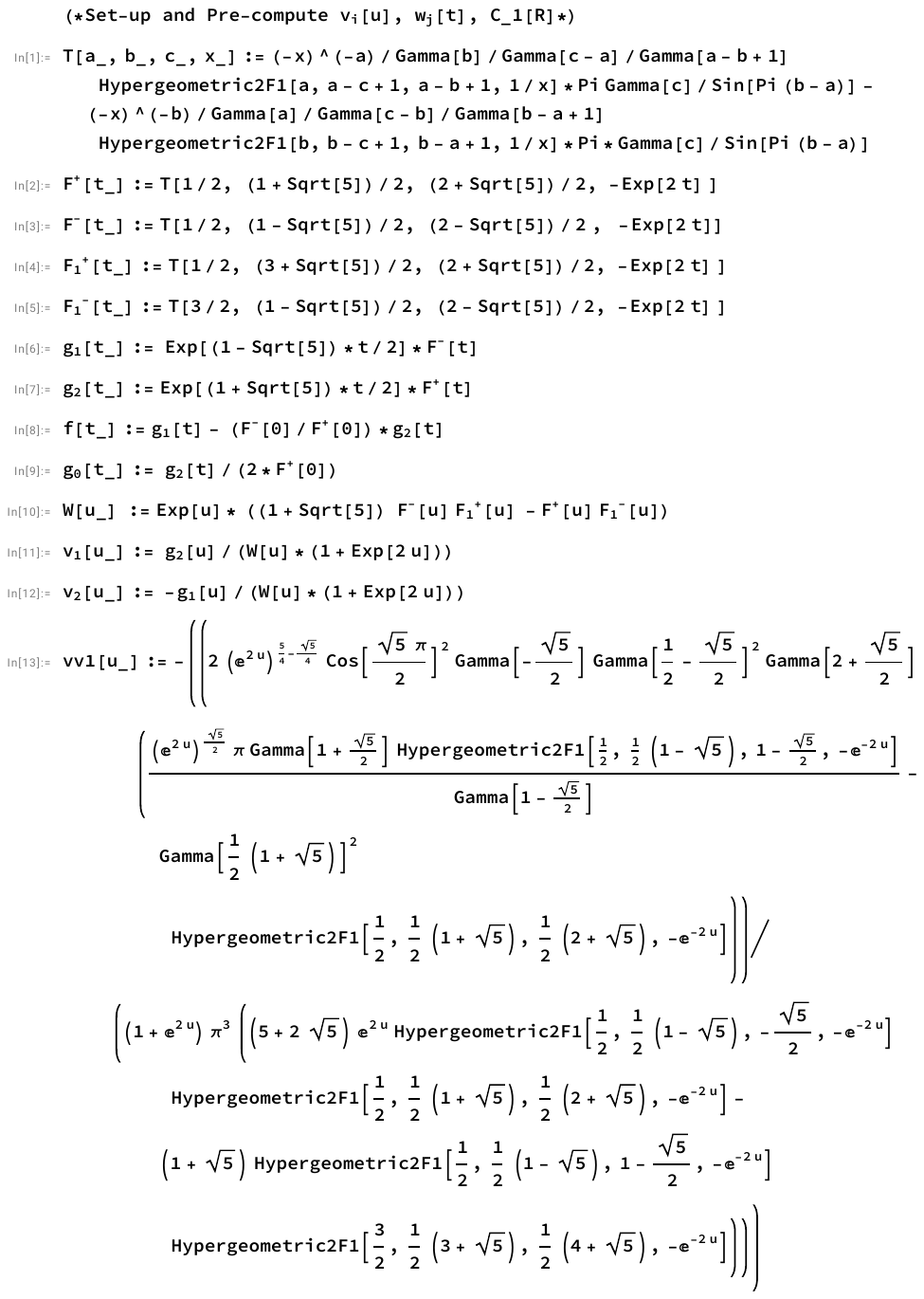}
\end{center}

\begin{center}
    \includegraphics[page=2, width=0.9\textwidth, trim=1.8cm 2cm 2cm 2.5cm,  clip]{Pre_compute.pdf}
\end{center}

\begin{center}
    \includegraphics[page=3, width=0.9\textwidth, trim=1.8cm 2cm 2cm 2.5cm,  clip]{Pre_compute.pdf}
\end{center}

\begin{center}
    \includegraphics[page=4, width=0.9\textwidth, trim=1.8cm 2cm 2cm 2.5cm,  clip]{Pre_compute.pdf}
\end{center}

\begin{center}
    \includegraphics[page=1, width=0.85\textwidth,trim=1.8cm 1.8cm 2cm 2.5cm, clip]{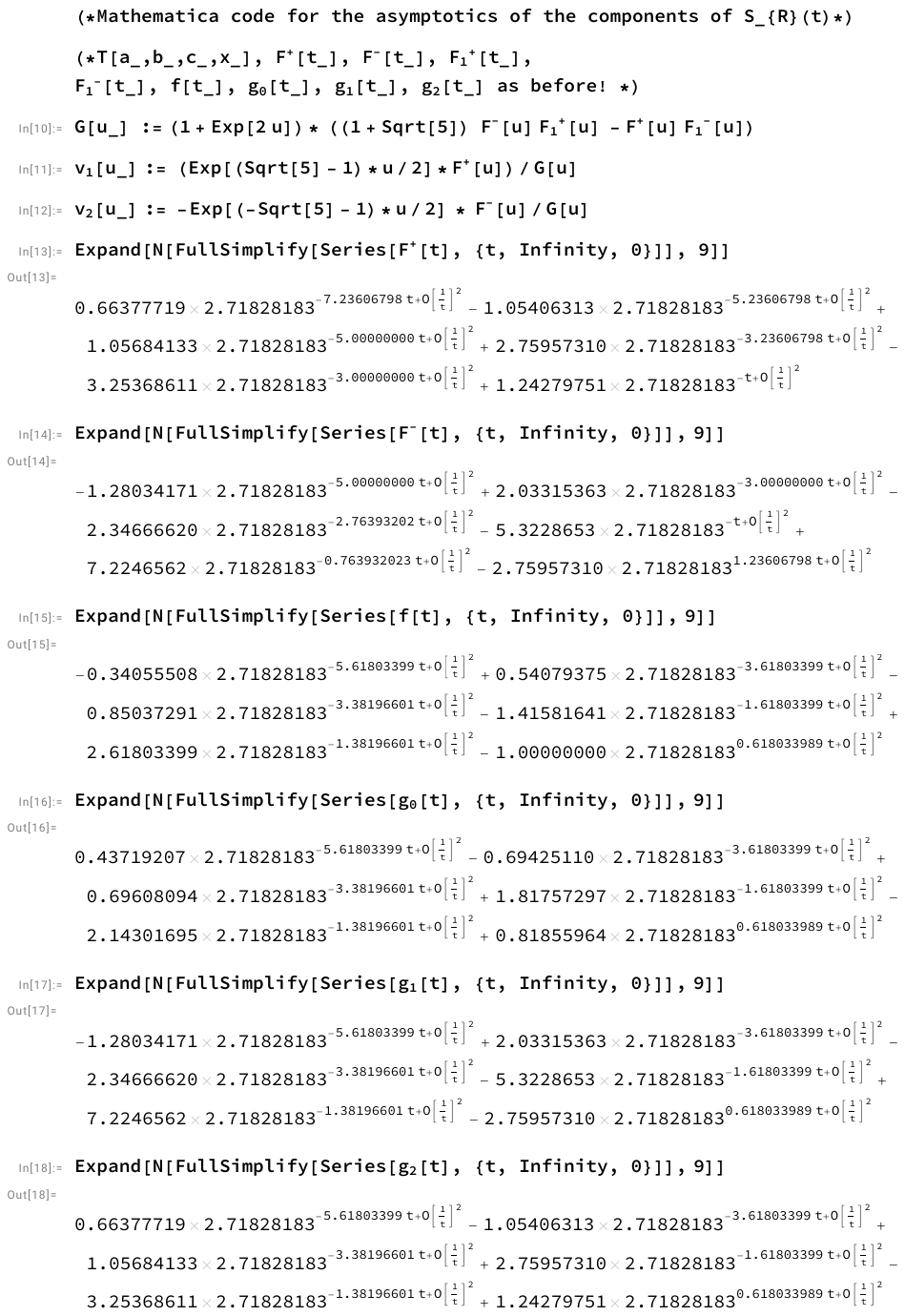}
\end{center} 

\begin{center}
    \includegraphics[page=2, width=0.85\textwidth,trim=1.8cm 1.8cm 2cm 2.5cm, clip]{30Jul25_asymp_display.pdf}
\end{center} 

\begin{center}
    \includegraphics[page=3, width=0.85\textwidth,trim=1.8cm 1.8cm 2cm 2.5cm, clip]{30Jul25_asymp_display.pdf}
\end{center} 

\begin{center}
    \includegraphics[page=4, width=0.85\textwidth,trim=1.8cm 13cm 2cm 2.5cm, clip]{30Jul25_asymp_display.pdf}
\end{center} 

\begin{center}
    \includegraphics[page=1, width=0.85\textwidth,trim=1.8cm 1.8cm 2cm 2.5cm, clip]{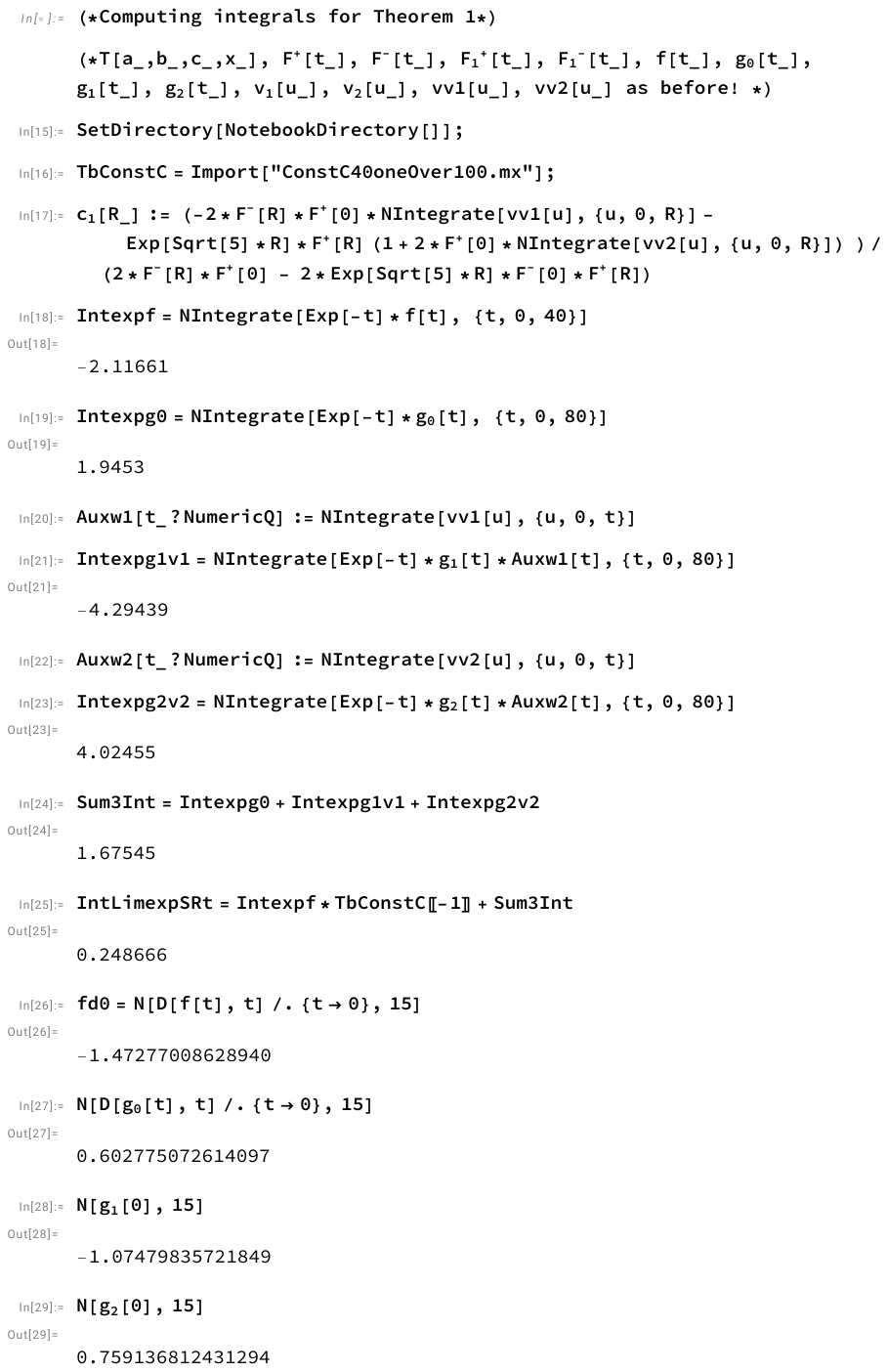}
\end{center} 

\begin{center}
    \includegraphics[page=2, width=0.85\textwidth,trim=1.8cm 1.8cm 2cm 2.5cm, clip]{30Jul_integral_useNInt_display_mode.pdf}
\end{center} 

\begin{center}
    \includegraphics[page=3, width=0.85\textwidth,trim=1.8cm 20cm 2cm 2.5cm, clip]{30Jul_integral_useNInt_display_mode.pdf}
\end{center}

     \end{document}